\documentclass[11pt,leqno]{article}
\usepackage[cp1251]{inputenc}
\usepackage[english,bulgarian]{babel}
\usepackage{amssymb}
\usepackage{amsfonts}
\usepackage{amsthm}
\usepackage{amsmath}

\newcommand{\X}{\mathfrak X}

\begin{document}
\setcounter{page} {144}
\begin{center} 
\bf {\footnotesize MATHEMATICS AND EDUCATION IN MATHEMATICS, 1990 \\
Proceedings of the Nineteenth Spring Conference of the Union of \\
Bulgarian Mathematicians } \\
{\footnotesize Sunny Beach, April 6--9, 1990} \\
{\footnotesize Sofia, BAS, 1990}
\end{center}
\vspace*{.51cm}
\begin{center}
\footnotesize 
 SUBCLASSES OF THE CONFORMAL ALMOST CONTACT METRIC MANIFOLDS
\end{center}
\begin{center}
Milen J. Hristov, Valentin A. Alexiev
\end{center}
\vspace*{.41cm}
\hspace*{1.5cm}
$\begin{array}{l} \mbox{\small A \, classification\, scheme\, of\, the\, conformal\, almost\, contact\, metric} \\ \mbox{\small manifolds\, with\, respect\, to\, the\, covariant\, derivative\, of\, the\, Lee} \\ \mbox{\small form \, is\,  given. \,The \, subclasses\, of\, one\, basic\, class\, and\, their } \\ \mbox{\small exact\, characterizations\, by\, the\, maximal\, subgroups\, of\, the\, contact} \\ \mbox{\small conformal\, group\, preserving\, itself\, are\, found.}
\end{array} $
\\ \\  \par 
  1. {\footnotesize PRELIMINARIES}. 
\small
Let $V$ be a $(2n+1)$-dimensional real vector space with almost contact metric structure $(\varphi ; \xi ;\eta ;g)$, where $\varphi$ is a tensor of type $(1,1),\, \xi$ is a vector, $\eta$ is a covector and $g$ is a positive definite metric such that
\[ 
\varphi ^2=-id+\eta\otimes\xi , \quad \varphi (\xi)=0, \quad g(\xi, \xi )=1,\quad g\circ \varphi =g-\eta\otimes\eta.
\]
\quad\quad The operators $h$ and $v$ on $V$ are defined by
\\ $ \hspace*{1.5cm}
h=-\varphi ^2, \quad v=\eta\otimes \xi .
$ \\
\quad\quad The decomposition $V=hV\oplus vV$ is orthogonal and invariant under the action of $U(n)\times 1$.
 \vskip 2mm
  2. {\footnotesize THE SPACE OF THE TENSORS OF TYPE (0,2) ON AN ALMOST CONTACT METRIC VECTOR SPACE}.  Let $\mbox{\footnotesize $L$}=V^*\otimes V^*$ denote the space of the tensors of type $(0,2)$ on $V^{2n+1}(\varphi , \xi ,\eta ,g)$. The metric $g$ induces a natural inner product in $\mbox{\footnotesize $L$} $. The standard representation of $U(n)\times 1$ in $V$ induces an associated representation of $U(n)\times 1$ in $\mbox{\footnotesize $L$} $. The operators $S$, $A$, $h$, $v$ and $w$ on $\mbox{\footnotesize $L$} $ are defined by
\vskip 2mm \qquad $
\begin{array}{l}
SL(x,y)= [L(x,y)+L(y,x) ]/2 , \\
AL(x,y)=[L(x,y)-L(y,x) ]/2 , \\
hL(x,y)=L(hx, hy), \\
vL(x,y)=\eta (x).L(\xi ,y)+ \eta (y).L(x, \xi )-2\eta (x)\eta (y).L(\xi ,\xi) ,\\
wL(x,y)=\eta (x)\eta (y).L(\xi ,\xi) ,
\end{array}
$ \\ \\
for arbitrary $L\in \mbox{\footnotesize $L$}$, $x$ and $y$ in $V$.
\par The traces $\alpha$ and $\beta$ of a trensor $L$ in $\mbox{\footnotesize $L$}$ are defined by 
\[
\alpha =g (e^i,e^j).hL(e_i, e_j), \quad \beta =g (e^i,e^j).hL(e_i, \varphi e_j),
\]
where $\{e_i \}, i=1,\ldots , 2n+1$ is a basis of $V$ and $\{e^i \}$ is the dual basis.
\vskip 2mm The tensors $L_i(L), \, i=1,\ldots ,9 $ associated with an arbitrary $L\!\!\in\!\! \mbox{\footnotesize $L$}\,$ are defined 
\[ 
\begin{array}{l}
L_1(L)(x,y)=\alpha.hg(x,y), \\ 
L_2(L)(x,y)=\frac{1}{2}[ ShL(x,y)+ ShL\circ \varphi (x,y) ]- L_1(L)(x,y), \\
L_3(L)(x,y)=\frac{1}{2}[ ShL(x,y)- ShL\circ \varphi (x,y) ] , \\
L_4(L)(x,y)=\beta.hg(x,\varphi y), \\
L_5(L)(x,y)=\frac{1}{2}[ AhL(x,y)+ AhL\circ \varphi (x,y) ]- L_4(L)(x,y), \\
L_6(L)(x,y)=\frac{1}{2}[ AhL(x,y)- AhL\circ \varphi (x,y) ] , \\
L_7(L)(x,y)=SvL(x,y), \\
L_8(L)(x,y)=AvL(x,y), \\
L_9(L)(x,y)=wL(x,y), 
\end{array}\leqno{(1)}
\]
for all $x,y$ in $V$.
\par The subspaces $\mbox{\footnotesize $L_i$}$ of $\mbox{\footnotesize $L$}$ are determined by the conditions
\vskip 2mm \qquad \quad $\mbox{\footnotesize $L_i$}=\{ L\in \mbox{\footnotesize $L$}: L=L_i(L)\}, i=1,\ldots ,9 $.
\vskip 2mm Since for any $L\in \mbox{\footnotesize $L$}$ are valid the following equations 
\vskip 2mm \qquad \quad $L=ShL+AhL+SvL+AvL+wL$,
\vskip 2mm \qquad \quad $ShL=L_1(L)+L_2(L)+L_3(L)$,
\vskip 2mm \qquad \quad $AhL=L_4(L)+L_5(L)+L_6(L)$,
\\ 
it is not difficult to prove the following
\vskip 2mm
 {\footnotesize THEOREM} 1. The decomposition 
 \[ \hspace*{-8cm}
 \mbox{\footnotesize $L$}=\mathop{\oplus}\limits_{i=1}^9 \mbox{\footnotesize $L_i$}
 \]
 is orthogonal and invariant under the action of $U(n)\times 1 $. The corresponding components of a tensor $L$ in 
 $\mbox{\footnotesize $L$}$ are $L_i(L), i=1,\ldots ,9$.
\vskip 2mm
  3. {\footnotesize APPLICATIONS TO ALMOST CONTACT METRIC MANIFOLDS}.  Let $M^{2n+1}(\varphi , \xi ,\eta ,g)$ be a $(2n+1)$-dimensional almost contact metric manifold where $\varphi$ is a tensor field of type (1,1), $\xi$ - a vector field, $\eta$ - 1-form and g a Riemannian metric such that
 \[
 \varphi ^2=-id + \eta\otimes \xi,\quad \varphi (\xi )=0,\quad \eta\circ\varphi=0,\quad g(\xi ,\xi)=1,\quad g\circ\varphi=g-\eta\otimes\eta .
 \]
The tangential space $T_pM$ at $M$ in $p\in M$ is an almost contact metric vector space.
\par Let $\nabla$ be the Levi-Civita connection on $M$. The fundamental 2-form $\Phi$ is given by $\Phi (x,y)=g(x,\varphi y)$ for all tangent vectors $x,y\in T_pM$, $p\in M$. We denote
\[
F=-\nabla\Phi .
\]
\quad The following 1-forms are associated with $F$:
\[
f(z)=\sum_{i=1}^{2n+1}F(e_i,e_i,z),\quad f^*(z)=\sum_{i=1}^{2n+1}F(e_i,\varphi e_i,z),\quad \omega (z)=F(\xi,\xi,z),
\]
where $\{e_i \}$, $i=1,\ldots , 2n+1$ is an arbitrary orthonormal basis of $T_pM$, $z\in T_pM$ and $p\in M$. We denote $\Omega=\X ^*M$ -- the space of 1-forms on $M$ and $L_2^0M$ -- the space of tensor fields of type (0,2) on $M$. The covariant derivative of a form defines the map
\[
L:\theta\in\Omega \longrightarrow L(\theta )=\nabla\theta \in L_2^0M.
\]
\quad Thus the subspaces $\mbox{\footnotesize $L$}_i$ in the decomposition of the space $\mbox{\footnotesize $L$}$ given in Theorem\,1 induce the corresponding subspaces $\Omega _i$ of the space $\Omega$. More precisely
\[
\Omega_i=\{ \theta\in\Omega:L(\theta )=L_i(\theta )\}, \leqno{(2)}
\]
where $L_i(\theta )$ are the components of $L(\theta )$ in $\mbox{\footnotesize $L$}_i$ , i.e. $L_i(\theta )_p=L_i(L(\theta )_p)$. The subspaces corresponding to $\mbox{\footnotesize $L$}_i\oplus \mbox{\footnotesize $L$}_j$ will be denoted by $\Omega _i\oplus \Omega _j$, $i,j=1,\ldots, 9.$
\vskip 2mm The following proposition is well known.
\par  {\footnotesize PROPOSITION} 1. Let $\theta\in\Omega$. Then
\par i) $L(\theta)$ is a symmetric tensor field iff $\theta$ is closed;
\par ii) $L(\theta)$ is an anisymmetric tensor field iff the dual vector field corresponding to $\theta$ is a Killing vector field.
\par From   Theorem 1. and  Proposition 1. it follows
\vskip 2mm  {\footnotesize PROPOSITION} 2. Let $\theta \in \Omega$. Then
\par i) $\theta $ is closed iff $\theta\in \Omega_1\oplus\Omega_2\oplus \Omega_3\oplus \Omega_7\oplus \Omega_9$;
\par i) $\theta ^*$ is Killing vector field iff $\theta\in \Omega_4\oplus\Omega_5\oplus \Omega_6\oplus \Omega_8$,
\\ where $\theta ^*$ is the dual vector corresponding to $\theta$.
\vskip 3mm
  4. {\footnotesize CLASSIFICATION SCHEME FOR THE CLASS $W_1$ OF ALMOST CONTACT METRIC MANIFOLDS}.  \,\, Using the decomposition of the space of tensors having the same symmetries as the covariant derivative of the fundamental 2-form in [2] is given a classification scheme for the almostcontact metric manifolds containing 12 basic classes $W_i$, $i=1, \ldots ,12$. The manifolds in the class $W_1\oplus W_2\oplus W_3\oplus W_9$ are said to be coformal classes [4]. These manifolds are generated by Sasakian and cosymplectic manifolds by means of subgroups of the contact conformal group and they are the contact analogues of the conformally Keahler manifolds in the Hermitian geometry [3], [5].
 \par Analogously to the Hermitian case [1] by making use of (2) we get a classification scheme for conformal manifolds with respect to the covariant derivative of the Lee form. 
 \par The 1-form $\theta$ on a manifold $M^{2n+1}(\varphi, \xi ,\eta, g)$ defined by 
 \[
 \theta=\frac{f^*(\xi)}{2n}.\eta+\frac{1}{2(n-1)}.f\circ \varphi
 \]
 is called Lee form on $M$.
 \par We call a manifold in the class $W_i$ is in the subclass $W_{ij}$, $i=1,2,3,9$,\\
  $j=1,\ldots , 9$ if $\theta $ is in the subspace $\Omega _j$. The class corresponding to the subspace $\Omega_j\oplus \Omega _k$ will be denoted by $W_{ij}\oplus W_{ik}$.
 \par We consider the class $W_1$. The defining condition for this class is \\ 
 $F=\eta\otimes (\eta\wedge \omega)$. The Lee form $\theta $ on a manifold in $W_1$ is given by $\theta =\omega\circ\varphi$. In [5] is proved that $\theta$ is contact closed, i.e. $hd\theta =0$. Then Proposition 1 implies
 \[
 L_4(\theta)=L_5(\theta)=L_6(\theta)=0. \leqno{(3)}
 \]
 The subclass $W_1^0$ of $W_1$ consists of all manifolds with $\theta$ closed, i.e. $d\theta =0$. From Proposition 1 we have
 \[
 L_4(\theta)=L_5(\theta)=L_6(\theta)=L_8(\theta)=0. \leqno{(4)}
 \]
 Theorem 1 and equations (3), (4) imply 
 \par   {\footnotesize THEOREM} 2. The class $W_1$ and its subclass have the following subclasses:
 \[
 \begin{array}{l}
 W_1=W_{11}\oplus W_{12}\oplus W_{13}\oplus W_{17}\oplus W_{18}\oplus W_{19}, \\ {} \\
 W_1^0=W_{11}\oplus W_{12}\oplus W_{13}\oplus W_{17}\oplus W_{19}.
 \end{array}
 \]
\quad  5. {\footnotesize THE SUBCLASSES $W_{1j}$ OF THE CLASS $W_1^0$ AND THEIR EXACT CHARACTERIZATIONS}. 
 \\ From theorems 1,2 we obtain the defining conditions for the subclasses 
 \\ $W_{1j}=\{ M(\varphi, \xi ,\eta, g)\in W_1^0 : L(\theta ) =L_j(\theta ) \} , j=1,\ldots ,9$.
 \par A contact conformal transformation of the structure $(\varphi, \xi ,\eta, g)$ on an almost contact metric manifold $M$ is defined (see [3]) by 
 \[
(5)\quad c(u,v): M^{2n+1}(\varphi, \xi ,\eta, g)\longrightarrow M^{2n+1}(\overline{\varphi }=\varphi,\,\, \overline{\xi }={\rm e}^{-v}\xi \,\, ,\overline{\eta }={\rm e}^{v}\eta,\,\,\overline{g}={\rm e}^{2u} hg +{\rm e}^{2v}\eta\otimes\eta ) ,
\]
 where $u,\, v$ are differentiable functions on $M$. The set of all these transformations forms the contact conformal group $G$.
 \par Let $M(\varphi, \xi ,\eta, g)$ and $M(\overline{\varphi}, \overline{\xi} ,\overline{\eta}, \overline{g})$ be conformally related as in (5). The Levi-Civita connections $\nabla$ and $\overline{\nabla}$ of the both structures are related (see [4]) by
\[ (6) \,\, \begin{array}{rl}
2\overline{g}(\overline{\nabla}_x y, z)&=2{\rm e}^{2u}g(\nabla_x y, z)+2{\rm e}^{2u}\{ du(x)g(y,z)+du(y)g(x,z)-du(z)g(x,y)\} + \\ {} \\
{} & + \,\, 2\{[ {\rm e}^{2v}dv(x) - {\rm e}^{2u} du(x) ]\eta(y)\eta(z)+ [ {\rm e}^{2v}dv(y) - {\rm e}^{2u} du(y) ]\eta(x)\eta(z) - \\ {} \\
{} & - \,\, [ {\rm e}^{2v}dv(z) - {\rm e}^{2u} du(z) ]\eta(x)\eta(y) \} +({\rm e}^{2v}-{\rm e}^{2u}).[ 2\eta (\nabla_x y)\eta(z)- \\ {} \\ 
{} & - \,\, \eta(z).F(x,\xi,\varphi y) - \eta(y).F(x,\xi,\varphi z) - \eta(x).F(y,\xi,\varphi z) - \\ {} \\
{} & - \,\, \eta(z).F(y,\xi,\varphi x) + \eta(y).F(z,\xi,\varphi x) + \eta(x).F(z,\xi,\varphi y)] .
 \end{array} 
 \]
 \quad We have the following 
\par  {\footnotesize THEOREM} ([3],[5]). The maximal subgroup of $G$ -- $G_1$ preserving the class $W_1$ 
consists of the transformations of type (5) satisfying the condition $du=0$. The  
maximal subgroup $G_1^0$ of $G$ preserving the subclass $W_1^0$ of $W_1$ consists of the trans-\\ formations of $G_1$ satisfying the condition $dv(\xi)\! =\! 0$.
\par It follows from (6)
\vskip 2mm
 {\footnotesize LEMMA}. Let $M(\varphi, \xi ,\eta, g)$ and $M(\overline{\varphi}, \overline{\xi} ,\overline{\eta}, \overline{g})$ be manifolds in the class $W_1 $ and the structures $(\varphi, \xi ,\eta, g)$, 
$(\overline{\varphi}, \overline{\xi} ,\overline{\eta}, \overline{g})$ on $M$ are conformally related as in (5) with a \\ transformation $c(u,v)\in G_1$. Then the corresponding to both structures \\ $\{ \nabla ,\, L(\theta),\, L_i(\theta) \}$ and $\{ \overline{\nabla },\, \overline{L}(\overline{\theta }),\, \overline{L}_i(\overline{\theta}) \}$, $i=1,2,3,7,9$ are related by 
\[\qquad
\overline{\nabla}_xy=\nabla_xy - {\rm e}^{2v-2u}\eta(x)\eta(y)h({\rm grad}\, v)+[dv(x)\eta(y)+dv(y)\eta(x)-dv(\xi)\eta(x)\eta(y) ]\xi\,\, ;
\]
\[ \hspace{-3cm}
\overline{L}(\overline{\theta })=L(\theta)-L(dv\circ \varphi)+{\rm e}^{2v-2u}.\theta ({\rm grad}\, v).\eta\otimes\eta ; \leqno{(7)}
\]
\[
\overline{L}_i(\overline{\theta })=L_i(\theta)-L_i(dv\circ \varphi)+{\rm e}^{2v-2u}.\theta ({\rm grad}\, v).\eta\otimes\eta ,\quad i=1,2,3,7,9 \leqno{(8)}
\]
for arbitrary $x,y$ in $T_pM$, $p\in M$.
\par It is not difficult to verify that the sets
\[
G_{1i}^0=\{ c(u,v)\in G: du=0, dv (\xi)=0, L(dv\circ\varphi )=L_i(dv\circ\varphi ) \},\quad i=1,2,3,9
\]
are subgroups of the contact conformal group.
 \par  {\footnotesize THEOREM} 3. The maximal subgroups of the contact conformal group preserving the subclasses $W_{1i}$ of $W_1^0$ are $G_{1i}^0,\quad (i=1,2,3,9)$.
 \par  {\footnotesize PROOF}.  Let $M(\varphi, \xi ,\eta, g)$ be in $W_{1i}$, $c(u,v)\in G_{1i}^0$ and $(\overline{\varphi}, \overline{\xi} ,\overline{\eta}, \overline{g})=c(u,v)(\varphi, \xi ,\eta, g)$. It is known that $M(\overline{\varphi}, \overline{\xi} ,\overline{\eta}, \overline{g})$ is in $W_i^0$. It follows from (7), (8) and the defining conditions of $G_{1i}^0$ that $\overline{L}(\overline{\theta})=\overline{L}_i(\overline{\theta})$, i.e. $M(\overline{\varphi}, \overline{\xi} ,\overline{\eta}, \overline{g})$ is in $W_{1i}$.
 \par For the inverse let $M(\varphi, \xi ,\eta, g)$ and $M(\overline{\varphi}, \overline{\xi} ,\overline{\eta}, \overline{g})$ be two manifolds in $W_{1i}$ which are conformally related by a trasformation $c(u,v)\in G_1^0$.
 Then (7) and (8) imply $L(dv\circ \varphi)=L_i(dv\circ \varphi)$, i.e. $c(u,v)\in G_{1i}^0$, $i=1,2,3,7,9$.
 \par In [5] is prooved that an almost contact metric manifold $M$ is in $W_1^0$ iff $M$ is contact conformally related to a cosymplectic manifold by a transformation of the subgroup $G_1^0$. We get an analogously characterization for the classes $G_{1i}^0$ and the class of cosymplectic manifolds.
 \par  {\footnotesize THEOREM} 4. An almost contact metric manifold is in the class $W_{1i}$ iff the structure of the manifold is locally conformal to a cosymplectic structure by a transformation of the group $G_{1i}^0$, $i=1,2,3,7,9$.
 \par  {\footnotesize PROOF}.  If $M(\varphi, \xi ,\eta, g)$ is a cosymplectic in [3] is proved that by an arbitrary transformation of $G_1^0$ we obtain $M(\overline{\varphi}, \overline{\xi} ,\overline{\eta}, \overline{g})\in W_1^0$. Then (7) and (8) imply \\ 
 $M(\overline{\varphi}, \overline{\xi} ,\overline{\eta}, \overline{g})\in W_{1i}$ if $c(u,v)\in G_{1i}^0$.
 \par For the inverse let $M(\varphi, \xi ,\eta, g)$ be in $W_{1i}$. Since the Lee form  $\theta=\omega\circ\varphi$ is closed, the Poincare's lemma implies that there exists locally function $v$ such that $\theta = dv$. In [3] is proved that by a transformation $c(u,v)$ with the above function $v$ and function $u:du=0$ we obtain $M(\overline{\varphi}, \overline{\xi} ,\overline{\eta}, \overline{g})$ is a cosymplectic manifold, i.e. $\overline{F}=0$ and hence $\overline{\theta}=0$. Thus (7) and (8) imply  $L(dv\circ\varphi)=L_i(dv\circ\varphi)$ and hence the transformation $c(u,v)$ is in $G_{1i}^0$.
 \begin{center}
{\small R E F E R E N C E S }
 \end{center}
 1. F.Tricherri, I. Vaisman. On some 2-dimensional Hermitian manifolds. Math. Z., \\ \hspace*{3mm} \underline{192} (1986), 205-216.
 \\
 2. V.Alexiev, G.Ganchev. On the classification of the almost contact metric ma-
 \\ \hspace*{3mm} nifolds. Mathematics and Educ. in Math., Sunny Beach 1986, 155-161.
 \\
 3. V.Alexiev, G.Ganchev. Contact conformal transformations of the almost con-
 \\ \hspace*{3mm} tact metric structures. Comp. rend. Acad.bulg. Sci., \underline{39} (1986), 27-30.
 \\
 4. V.Alexiev, G.Ganchev. Canonical connection on a conformal contact metric ma-
 \\ \hspace*{3mm} nifold. Annuaire de l'Universite de Sofia, \underline{81}, Fac. Math. et Inf. (1987) in 
 \\ \hspace*{3mm} press.
 \\
 5. V.Alexiev, G.Ganchev. On a class of the almost contact metric manifolds conformal-
 \\ \hspace*{3mm} ly related to cosymplectic manifolds. Comp. rend. Acad.bulg. Sci., \underline{41} (1988), 
 \\ \hspace*{3mm} 21-24.
\end{document}